\documentclass[11pt, reqno]{amsart}
\usepackage{amsthm,amsmath,amsfonts,amssymb,color}
\usepackage[bookmarks]{hyperref}
\reversemarginpar

\newtheorem{thm}{Theorem}[section]

\newtheorem{cor}[thm]{Corollary}

\newtheorem{preremark}[thm]{Remark}
\newenvironment{remark}{\begin{preremark}\rm}{\medskip \end{preremark}}

\numberwithin{equation}{section}

\newcommand{\R}{\mathbb R}

\DeclareMathOperator{\Vol}{Vol}

\newcommand{\grad} {\nabla}

\newcommand{\dd} {\mathrm{d}}

\DeclareMathOperator{\dv}{div}

\DeclareMathOperator{\Ric}{Ric}

\def\H{\mathbb H^{2}(-a^{2})}
\def\be{\begin{equation}}
\def\ee{\end{equation}}

\begin{document}
\title[Separation Equation]{An ODE for  boundary layer separation 
on a sphere and a hyperbolic space}

\author[Chan]{Chi Hin Chan}
\address{Department of Applied Mathematics, National Chiao Tung University,1001 Ta Hsueh Road, Hsinchu, Taiwan 30010, ROC}
\email{cchan@math.nctu.edu.tw}
\author[Czubak]{Magdalena Czubak}
\address{Department of Mathematical Sciences, Binghamton University (SUNY),
Binghamton, NY 13902-6000, USA}
\email{czubak@math.binghamton.edu}
\author[Yoneda]{Tsuyoshi Yoneda}
\address{Department of Mathematics, Hokkaido University,
Kita 10, Nishi 8, Kita-Ku, Sapporo, Hokkaido, 060-0810, Japan}
\email{yoneda@math.sci.hokudai.ac.jp}

\begin{abstract}
Ma and Wang derived an equation linking the separation location and times for the boundary layer separation of incompressible fluid flows. The equation gave a necessary condition for the separation (bifurcation) point.  The purpose of this paper is to generalize the equation to other geometries, and to phrase it as a simple ODE. Moreover we consider the Navier-Stokes equation with the Coriolis effect, which is related to the presence of trade winds on Earth.
\end{abstract}
\date{\today}
\maketitle


 \section{Introduction}
In the beginning of 20th century,  Prandtl proposed the boundary layer theory.  Since then there has been a lot of extensive developments in the theory (see Rosenhead \cite{R} 
 for example).
In general, the laminar flow in the boundary layer 
should be governed by a boundary layer equation, which is deduced from the Navier-Stokes equations.
The existence of the singularity in the steady boundary layer flow along fixed wall has led to important advances in the understanding 
of the steady boundary layer separation. In this point of view, 
Van Dommelen and Shen \cite{DS} made a key observation of shock singularities with numerical computations.  In the beginning of  21th century, Ghil, Ma and Wang \cite{GMW, GMW2, M, MW} have developed a rigorous theory on the boundary layer separation of incompressible fluid flows. 
Their articles are oriented towards the structural bifurcation and boundary layer separation of the solution to the Navier-Stokes equations.
In particular, in \cite{MW} authors established a simple equation, which they call a ``separation equation," linking the separation location and times. Furthermore, they showed that the structural bifurcation occurs at a degenerate singular point with integer index of the velocity field at the critical bifurcation time. Their theory is based on the classification of the detailed orbit structure of the velocity 
field near the bifurcation time and location (see also \cite{B}).
 On the other hand, Ghil, Liu, Wang and Wang \cite{GLWW}  gave a new rigorous argument of ``adverse pressure gradient" mathematically under certain conditions.  The conditions were consistent with the careful numerical experiment also found in \cite{GLWW}.  The appearance of the adverse pressure gradient is well known to be the main mechanism for the boundary-layer separation in physics. 
 
The purpose of this paper is to obtain the separation equation of Ma and  Wang's in other geometries and to phrase it as an ODE.  In order to state our main result, we need to explain Ma and Wang's separation equation precisely.
Let $K$ be a compact domain in $\R^2$ with $C^{r+1}$ boundary, $\partial K$, for $r\geq 2$.  Consider the Navier-Stokes equation on $K$ given by

 \be\label{NS}
 \begin{split}
u_t+\nabla_u u-\Delta u+\grad p&=0,\\
\dv u&=0,\\
u_{|\partial K}&=0,\\
u(x,0)&=\phi(x), \quad \phi_{|\partial K}=0.
\end{split}
 \ee
Since we only consider the flow near the boundary, we can replace $K$ by $\mathbb{R}^2-K$.

We call a point $p\in \partial K$ ``$\partial$-regular point of $u$" if the normal derivative of the tangential component of $u$ at $p$ is nonzero, i.e., $\partial (u\cdot \tau)(p)/\partial n\not=0$, otherwise, 
$p\in\partial K$ is called a $\partial$-singular point (bifurcation point) of $u$.

\begin{thm}\cite{MW}\label{MWthm}   
Let $K$ be a compact domain in $\R^2$ with $C^{r+1}$ boundary, $\partial K$, for $r\geq 2$.  Let $p_0 \in \partial K$, and $t_0\geq 0$.  If $(p_0,t_0)$ is a $\partial$-singular point  (bifurcation point) of the solution $u$ of \eqref{NS}, then 
\be\label{mw}
\frac{\partial \phi_\tau(p_0)}{\partial n}=\int^{t_0}_0\nabla \times \Delta u-k\Delta u\cdot \tau \dd t,
\ee
where $\nabla \times \Delta u=\partial_\tau (\Delta \cdot n)-\partial_n(\Delta u\cdot \tau)$, and $k(p_0)$ is the curvature of $\partial K$ at $p_0$.
\end{thm}
The equation \eqref{mw} is called the ``separation equation."
  
We now give an example which tells us imposing inflow profile is useful. 
A wind turbine system consisting
of a diffuser shroud with a broad-ring at the exit periphery and a wind turbine inside it was developed by Ohya and Karasudani \cite{OK}.
Their experiments show that a diffuser-shaped (not nozzle-shaped) structure  can accelerate the wind at the entrance of the body.  This is called a ``wind-lens phenomena."  A strong vortex formation with  a low-pressure region is created behind the broad brim.   The wind flows into a low-pressure region, and the wind velocity is increased more near the entrance of the diffuser.  In general, creation of a vortex needs  separation phenomena near a boundary (namely, bifurcation phenomena), and 
 before separating from the boundary, the flow moves toward reverse direction near the
boundary against the laminar flow (inflow) direction.
In order to consider such phenomena in pure mathematics, imposing inflow profile at the entrance should be reasonable.

We moreover consider the situation on a sphere and a hyperbolic space (we can easily deduce the ODE in the Euclidean case).  In the case of the sphere, one motivation comes from studying the flow on the Earth (see Corollary \ref{inflow_cor}).

Now, we write the equation on a Riemannian manifold, $M$, where $M$ is taken either to be a sphere $S^2(a^2)$ or a hyperbolic space $M=\H$.  We write the equation in the language of differential $1$-forms as follows.

Let $O$ be the base point in $M$.  Let $(r,\theta)$ be the normal polar coordinates on $M$. Then we have the following orthonormal moving frame  
\begin{align}
e_1&=\partial_r,\\
e_2&=\frac 1{s_a(r)} \partial_\theta,
\end{align}
where
$s_a(r)=\frac{\sin(ar)}{a}$ if $M=S^2(a^2)$ or  $s_a(r)=\frac{\sinh(ar)}{a}$ if $M=\H$.  
 We also introduce $c_a(r)$, where $c_a(r)=\cos(ar)$ if $M=S^2(a^2)$ or  $c_a(r)={\cosh(ar)}$ if $M=\H$.  Note
\be\label{ds}
\frac{\dd}{\dd r}s_a(r)=c_a(r).
\ee
For simplicity, in the sequel, we omit the writing of subscripts $a$ in $s_a$ and $c_a$ and simply write $s, c$.
The associated dual frame to $\{e_1,e_2\}$ can be written as
\begin{align}
e^1&=\dd r,\\
e^2&=s(r)\dd \theta.
\end{align}
Hence the volume form on $M$ is given by  $\Vol_M=e^1\wedge e^2=s\dd r\wedge \dd \theta$.
Let $\nabla$ be the Levi-Civita connection on $M$.  We have
\begin{align}
\nabla_{\partial_r}\partial_r&=0,\\
\nabla_{\partial_r}\partial_\theta&=\nabla_{\partial_\theta}\partial_r= c e_2=\frac{c}{s}\partial_\theta,\\
\nabla_{\partial_\theta}\partial_\theta&=-cs\partial_r.
\end{align}
These imply
\be\label{table2}
\nabla_{e_1}e_1=\nabla_{e_1}e_2=0,\quad\nabla_{e_2}e_2=-\frac cs e_1,\quad \nabla_{e_2}e_1=\frac cs e_2.
\ee
 Let $d$ be the distance function on $M$.  Define and obstacle $K$ on $M$ by $K=\overline{(B_O(\delta))}=\{p\in M: d(p,O)\leq \delta \}$.  
Consider a smooth vector field $u$ defined on a neighborhood near $\partial K$.  Then $u$ can be written as
 \[
 u=u_r e_1+u_\theta e_2,
 \] 
 for some locally defined smooth functions $u_r, u_\theta$.  By ``lowering the index" we can obtain a 1-form $u^\ast=u_r e^1+u_\theta e^2$.  For simplicity we just write $u$ for both the vector field and the 1-form.
Recall the Hodge star operator, $\ast$, is a linear operator that sends $k$-forms to $n-k$-forms and is defined by  
\be\label{star}
\begin{split}
\alpha \wedge \ast \beta &=g(\alpha ,\beta)\Vol_M.
\end{split}
\ee
Then
\be\label{star2}
\ast \ast \alpha =(-1)^{nk+k}\alpha,
\ee
where $n$ is the dimension of the manifold, and $k$ the degree of $\alpha$.
Here, by a direct computation
\[
\ast e^1=e^2,\quad \ast e^2=-e^1,\quad  \ast \Vol_M=1.
\]
Recall
\be\label{star3}
\dd ^\ast \alpha=(-1)^{nk+n+1}\ast \dd \ast\alpha.
\ee
So for two dimensional manifolds we have $\dd^\ast=-\ast \dd\ast$ . 
Then the Navier-Stokes equation on $M-K$ is given by
   \be\label{NSM}
 \begin{split}
u_t+\nabla_u u-\Delta u-2\Ric u+\dd p&=0,\\
\dd^\ast u&=0,\\
u_{|\partial K}&=0,\\
u(x,0)&=u_0(x), \quad {u_0}_{|\partial K}=0.
\end{split}
 \ee
For fixed $p_0\in\partial K$, and $u$ a solution of $\eqref{NSM}$,
let us give  the key parameters
\begin{eqnarray*}
k=k_{a,\delta}&:=&\frac{c_a(\delta)}{s_a(\delta)},\\
\alpha_1(t)&:=&\partial_ru_\theta(t,p_0),\\
\alpha_2(t)&:=&\partial_r^2u_\theta(t,p_0),\\
\alpha_3(t)&:=&\partial_r^3u_\theta(t,p_0),\\
\eta(t)&:=&\frac{1}{s^2(\delta)}\partial_r\partial_\theta^2u_\theta(t,p_0).
\end{eqnarray*}
Note that $k$ includes both curvature of the manifold and curvature of the boundary.
Our main theorem is the following:
\begin{thm}\label{mainthm}
Let $\alpha_1(0)>0$ (initial data), and $\alpha_2(t)$, $\alpha_3(t)$ and $\eta(t)$ be given functions. Then $\alpha_1(t)$ satisfies the following ODE:
\begin{equation*}
\partial_t\alpha_1(t)=-k^2\alpha_1(t)+\alpha_3(t)+2k\alpha_2(t)+2\eta(t).
\end{equation*}
\end{thm}
\begin{remark}\ We give five remarks.

\begin{itemize}

\item

A $\partial$-singular point (bifurcation point) occurs at $t_0$ iff a function $\alpha_1(t)$ satisfies $\alpha_1(t_0)=0$. 

\item
The above result is a generalization of  \cite{Y} which is considered in the Euclidean space $\mathbb{R}^2$.
\item

We can regard $\alpha_1(t)$, $\alpha_2(t)$ and  $\alpha_3(t)$ as a part of the inflow profile. However 
$\eta(t)$ is not. Let us be more precise. Choose $\tilde p\in\partial K$ close to $p_0\in\partial K$,
and let
\begin{equation*}
\tilde K:=\{p\in M-K : d(p,\tilde p)<d(p_0,\tilde p)\}.
\end{equation*}
Then $\alpha_1(t)$, $\alpha_2(t)$ and $\alpha_3(t)$ can be determined by $u(\cdot,t)$ on $\partial\tilde K\cap K$ near $p_0\in\partial K$ (boundary value).
$\eta(t)$ can be determined by $u(\cdot,t)$ in $\tilde K\cap K$ near $p_0\in\partial K$ (interior flow).

\item
 
We can find a geometric meaning of $\eta(t)$ (see also \cite{Y}).
 
\begin{itemize}

\item
Convexing streamlines:\ 
We can see (geometrically) convexing streamlines near the boundary iff 
\begin{equation*}
\eta(t)
< 0.
\end{equation*}

\item
Almost parallel streamlines:\ 
We can see (geometrically) almost parallel streamlines near the boundary iff 
\begin{equation*}
\eta(t)
= 0.
\end{equation*}

\item

Concaving streamlines:\ 
We can see (geometrically) concaving streamlines near the bounday iff 
\begin{equation*}
\eta(t)
> 0.
\end{equation*}

\end{itemize}

\item

It is reasonable to assume  $u_\theta$ does not grow polynomially for $r$ direction (this is due to the observation of  ``boundary layer," since the flow should be a uniform one  away from the boundary).  Thus, it should be  reasonable to focus on the following two cases:
\begin{itemize}

\item (Poiseuille type profile)\ 
$-k^2\alpha_1(t)+2k\alpha_2(t)<0$ ($\alpha_1(t)>0$, $\alpha_2(t)<0$) and $\alpha_3(t)$ is small comparing with $\alpha_1(t)$ and $\alpha_2(t)$.

\item (Before separation profile)\  $2k\alpha_2(t)+\alpha_3(t)<0$ ($\alpha_2(t)>0$, $\alpha_3(t)<0$)  and $\alpha_1(t)$ is small comparing with $\alpha_2(t)$ and $\alpha_3(t)$. 

\end{itemize}

In this  point of view, the well-known physical phenomena of ``adverse pressure gradient" occurs
in ``before separation profile," since $\alpha_2(t)>0$ and $\dd p=\Delta u$ on the boundary.
\end{itemize}
\end{remark}
Our method can be applied to geophysics, in particular, to the ``trade winds'' on Earth. The trade winds are the easterly surface winds that can be found in the tropics, within the lower portion of the Earth's atmosphere near the equator. The Coriolis effect is responsible for deflecting the surface air, which flows from subtropical high-pressure belts toward the Equator, toward the west in both hemispheres. In the corollary below, we consider the Navier-Stokes equation with the Coriolis effect on a rotating sphere.  The equation is
\begin{eqnarray}
u_t+\nabla_{u}u-\triangle u -2 \Ric u  +\beta \cos (ar) * u  + \dd p & =& 0,\label{ceq1} \\
\dd^* u & =& 0,\\
u(x,0)&=&u_0\label{ceq2},
\end{eqnarray}
where $\beta\in\mathbb{R}$ is a Coriolis parameter.  The term $\beta \cos (ar) * u$ represents the effect upon the velocity $u$ due to the rotation of the sphere with constant speed $\beta$.
It is worthwhile to mention that the existence and uniqueness of parallel
laminar flows satisfying the stationary version of the above system has
been considered in \cite{CY}.   
More details on the Coriolis effect on a sphere can be found in \cite{MA}, and in the vorticity formulation, for example, in \cite{YY}. 
\begin{cor}\label{inflow_cor} 
Let $u$ satisfy \eqref{ceq1}-\eqref{ceq2} and the following conditions
\[
u_\theta|_{\partial K}=0, \quad u_r|_{\partial K}=\lambda_0\in \R,\quad \partial_r u_r|_{\partial K}=0.
\]
Then
$\alpha_1(t)$ satisfies 
  \begin{align*}
\partial_t \alpha_1(t)=-k(k+\lambda_0)\alpha_1(t)+(2k-\lambda_0)\alpha_2(t)+\alpha_3(t)+2\eta(t)+
\lambda_0\beta(a\sin(a\delta)-k\cos(a\delta)).
\end{align*}
\end{cor}
This is proved in section \ref{Earth}.
\begin{remark}
If $\lambda_0(a\sin(a\delta)-k\cos(a\delta))$ is strictly positive, and $\beta$ is sufficiently large compared with $\alpha_1(0)>0$, $\lambda_0$, $k$, $\alpha_2(t)$, $\alpha_3(t)$ and $\eta(t)$,
then $\alpha_1(t)$ can never be zero. This expresses that  the south (or north) flow deflects towards the east (or west).
Moreover we can find an asymptotic behavior of $\alpha_1(t)$: 
\begin{equation*}
\lim_{t\to\infty}\left|\alpha_1(t)- \frac{(2k-\lambda_0)\tilde \alpha_2+\tilde \alpha_3+2\tilde \eta+
\lambda_0\beta(a\sin(a\delta)-k\cos(a\delta))} {k(k+\lambda_0)}\right|=0,
\end{equation*}
if $\alpha_2(t)\to \tilde\alpha_2$, $\alpha_3(t)\to\tilde \alpha_3$ and  $\eta(t)\to\tilde \eta$.
\end{remark}

\section{Proof of the main theorem.}
We first prepare the necessary computations and then put them together in Subsection \ref{finalproof}.
\subsection{Divergence free condition in coordinates}
If $u$ is divergence free, then $\dd^\ast u=0$.  Compute
\begin{align*}
0=\dd^\ast u= -\ast \dd \ast u &=-\ast \dd \ast u\\
&=-\ast \dd (u_r e^2-u_\theta e^1)\\
&=-\ast (\partial_r(su_r)+\partial_\theta u_\theta) \dd r\wedge \dd \theta\\
&=-\frac 1s  (\partial_r(su_r)+\partial_\theta u_\theta).
\end{align*}
This implies
\be\label{df}
\partial_r(su_r)+\partial_\theta u_\theta=0.
\ee
In addition on $\partial K$, thanks to the no-slip boundary condition,  from \eqref{df} we can deduce 
\be\label{df1}
0=\partial_\theta u_\theta |_{\partial K}=\{-cu_r-s\partial_ru_r\}|_{\partial K}=-s\partial_r u_r|_{\partial K}.
\ee

\subsection{Computing normal and tangential components of $\Delta u$.}
The goal is to compute $g(\Delta u,e^1)$ and $g(\Delta u, e^2)$.  First, $-\Delta u= \dd \dd^\ast u+\dd^\ast \dd u=\dd^\ast \dd u$.  Next 
\begin{align*}
\dd u=\{\partial_r(s u_\theta)-\partial_\theta u_r\}\dd r\wedge \dd\theta=\frac 1s\{\partial_r (su_\theta)-\partial_\theta u_r\}\Vol_M.
\end{align*}
It follows
\begin{align*}
\Delta u&=-\dd^\ast\dd u\\
&=\ast \dd \ast \frac 1s\{\partial_r(s u_\theta)-\partial_\theta u_r\}\Vol_M\\
&=\ast \dd  \frac1s\{(\partial_r (su_\theta)-\partial_\theta u_r\}\\
&=\ast  \{\partial_r( \frac 1s\partial_r (su_\theta))-\partial_r(\frac 1s\partial_\theta u_r)\}e^1+\frac 1s\ast\{ \partial_\theta \partial_r (su_\theta)-\partial^2_\theta u_r\}   \dd\theta\\
&=\partial_r\{ \frac 1s(\partial_r(s u_\theta)-\partial_\theta u_r)\}e^2-\frac {1}{s^2}\{ \partial_r (s\partial_\theta u_\theta)-\partial^2_\theta u_r\}   e^1.
\end{align*}
Then
\be\label{Deltae1}
g(\Delta u, e^1)=\frac {1}{s^2}\{ \partial^2_\theta u_r-\partial_r (s\partial_\theta u_\theta)\}=\frac {1}{s^2}\{ \partial^2_\theta u_r-c\partial_\theta u_\theta-s\partial_r\partial_\theta u_\theta\},
\ee
which can be rewritten using \eqref{df} as follows
\begin{align}
g(\Delta u, e^1)&=\frac {1}{s^2}\{ \partial^2_\theta u_r-c\partial_\theta u_\theta+s\partial_r^2(s u_r)\}\nonumber\\
&=\frac {1}{s^2}\partial^2_\theta u_r-\frac{c}{s^2}\partial_\theta u_\theta\mp a^2u_r+2\frac cs\partial_r u_r+\partial^2_r u_r\label{Delta1},
\end{align}
where $\mp$ depends on the choice of $M$.  Here, and in the sequel, the upper sign refers to the sphere and the lower sign to the hyperbolic plane.
Next
\begin{align}
g(\Delta u, e^2)&= \partial_r\{ \frac 1s(\partial_r(s u_\theta)-\partial_\theta u_r)\}\nonumber\\
&= \partial_r\{ \frac cs u_\theta+\partial_r u_\theta-\frac 1s\partial_\theta u_r\}\nonumber\\
&=-\frac{1}{s^2}u_\theta+\frac cs \partial_r u_\theta +\partial_r^2 u_\theta+\frac{c}{s^2}\partial_\theta u_r-\frac 1s \partial_r\partial_\theta u_r
 \label{Deltae2}
\end{align}
since $\partial_r(c/s)=-(1/s^2)$ and $\partial_r (1/s)=-(c/s^2)$.

\subsection{Computing  $\frac {1}{s}\partial_\theta g(\Delta u,e^1)$ on $\partial K$. } 
First observe that on $\partial K$, from the no-slip boundary condition, \eqref{df1} and \eqref{Delta1} we have
\[
g(\Delta u, e^1)|_{\partial K}=\partial^2_r u_r{|_{\partial K}}.
\]
Hence
\be\label{Deltae3}
 \frac {1}{s}\partial_\theta g(\Delta u,e^1)|_{\partial K}=\frac 1s\partial_\theta \partial^2_r u_r{|_{\partial K}}.
 \ee
We need this formula to estimate the pressure term on the boundary.
\subsection{Computing  $\partial_r g(\nabla_u u,e_2)$ on $\partial K$.}
First, by the properties of the connection and \eqref{table2}
\begin{align}
\nabla_uu&=\nabla_{u_re_1+u_\theta e_2}(u_re_1+u_\theta e_2)\nonumber\\
&=u_r\nabla_{e_1}(u_re_1+u_\theta e_2)+u_\theta\nabla_{e_2}(u_re_1+u_\theta e_2)\nonumber\\
&=u_r(\partial_r u_r e_1+\partial_r u_\theta e_2)+u_\theta(\frac 1s\partial_\theta u_r e_1+u_r\frac cs e_2+\frac 1s\partial_\theta u_\theta e_2-u_\theta\frac cs e_1)\label{convect}.
\end{align}
Then
\begin{align}\label{convection1}
g(\nabla_u u, e_2)=u_r\partial_ru_\theta+u_\theta u_r \frac cs+\frac 1s u_\theta\partial_\theta u_\theta.
\end{align}
Differentiating and evaluating on the boundary and using the no-slip boundary condition, we reduce \eqref{convection1} to
\[
\partial_r g(\nabla_u u, e_2)|_{\partial K}=\partial_r u_r \partial_r u_\theta  +\frac 1s \partial_r u_\theta \partial_\theta u_\theta.
\]
But then the divergence free condition \eqref{df} again with the no-slip boundary condition imply
\be\label{convection}
\partial_r g(\nabla_u u, e_2)|_{\partial K}=\{\partial_r u_r \partial_r u_\theta  -\frac 1s \partial_r u_\theta \partial_r(su_r)\}|_{\partial K}=0.
\ee
\subsection{Computing  $\partial_r g(\nabla p,e_2)$ on $\partial K$ }
First
\[
\dd p=\partial_r p \dd r + \partial_\theta p \dd\theta= \partial_r p e^1 + \frac 1s\partial_\theta p e^2.
\]
Hence
\[
\nabla p= \partial_r p e_1+  \frac 1s\partial_\theta p e_2,
\]
and
\be\label{p1}
\partial_r g(\nabla p, e_2)=\partial_r (\frac 1s \partial_\theta p)=-\frac c{s^2}\partial_\theta p+\frac 1s \partial_r\partial_\theta p=-\frac cs g(\nabla p, e_2) +\frac 1s \partial_r\partial_\theta p.
\ee
Equivalently, we can write \eqref{p1} as
\be\label{p2}
\partial_r g(\nabla p, e_2)=-\frac cs g(\nabla p, e_2) +\frac {1}{s}\partial_\theta g(\grad p, e_1).
\ee

\subsection{Proof of the formula}\label{finalproof}
We now follow the proof in \cite{MW}, but without assuming that $p_0$ is a bifurcation point.  

Let $p_0\in\partial K$, and $t_0>0$. Begin by writing
\be\label{e1}
{\partial_r}|_{p_0} g(u(t_0, \cdot),e_2)-{\partial_r}|_{p_0} g(u_0,e_2)=\int^{t_0}_{0}\frac{\dd}{\dd t}\{{\partial_r}|_{p_0} g(u(t, \cdot),e_2)\}\dd t.
\ee
From \eqref{NSM} it follows
\begin{align*}
\frac{\dd}{\dd t}\{{\partial_r} g(u(t, \cdot),e^2)\}&={\partial_r}g(\Delta u(t, \cdot),e^2)+{\partial_r} g(2\Ric (u(t, \cdot)),e^2)\\
&\quad-{\partial_r} g(\nabla_{u(t,\cdot)}{u(t,\cdot)},e_2)-{\partial_r} g(\dd p(t,\cdot),e^2).
\end{align*}
We can simplify by using $\Ric (u)=a^2 u$ if $M=S^2(a^2)$ and $\Ric(u)=-a^2 u$ if $M=\H$, and write $\Ric(u)=\pm a^2 u$.    Also on the boundary we can use \eqref{convection} to write
 \begin{align*}
\frac{\dd}{\dd t}\{{\partial_r} g(u(t, \cdot),e_2)\}&={\partial_r}g(\Delta  u(t, \cdot),e^2)\pm2a^2{\partial_r} g(u(t, \cdot),e_2) -{\partial_r} g(\grad p(t,\cdot),e_2)\\
&={\partial_r}g(\Delta  u(t, \cdot),e^2)\pm2a^2{\partial_r}u_\theta -{\partial_r} g(\grad p(t,\cdot),e_2)
\end{align*}
Next from \eqref{p2} we have
 \begin{align*}
\frac{\dd}{\dd t}\{{\partial_r} g(u(t, \cdot),e_2)\}={\partial_r}g(\Delta  u(t, \cdot),e^2)\pm2a^2{\partial_r}u_\theta +\frac cs g(\nabla p(t,\cdot), e_2) -\frac {1}{s}\partial_\theta g(\grad p(t,\cdot), e_1).
\end{align*}
Since on $\partial K$
\be\label{pressureb}
\dd p=\Delta u, 
\ee
going back to \eqref{e1} we obtain (compare this with \eqref{mw})
\begin{align}
&{\partial_r}|_{p_0} g(u_0,e_2)- {\partial_r}|_{p_0} g(u(t_0, \cdot),e_2)\nonumber\\
&=-\int^{t_0}_{0}   {\partial_r}g(\Delta  u(t, p_0),e^2)\pm2a^2{\partial_r}u_\theta +\frac cs g(\nabla p(t,p_0), e_2) -\frac {1}{s}\partial_\theta g(\grad p(t,p_0), e_1)     \dd t\nonumber\\
&=-\int^{t_0}_{0}   {\partial_r}g(\Delta  u(t, p_0),e^2)-\frac {1}{s}\partial_\theta g(\Delta u(t, p_0), e^1)  +\frac cs g(\Delta u(t,p_0), e^2)    \pm2a^2{\partial_r}u_\theta \dd t.\label{e2}
\end{align}
To obtain the necessary and sufficient condition, we write \eqref{e2} more explicitly as follows.  From \eqref{Deltae2}
\begin{align}
{\partial_r}g(\Delta  u(t, p_0),e^2)
&=\{\partial_r (-\frac{1}{s^2}u_\theta+\frac cs \partial_r u_\theta +\partial_r^2 u_\theta+\frac{c}{s^2}\partial_\theta u_r-\frac 1s \partial_r\partial_\theta u_r)\}|_{(t,p_0)}\nonumber\\
 &=\{-\frac 2{s^2}\partial_r u_\theta+\frac cs \partial^2_r u_\theta +\partial_r^3 u_\theta-\frac 1s \partial^2_r\partial_\theta u_r\}|_{(t,p_0)}\label{e3}.
\end{align}
And again from \eqref{Deltae2}
\be\label{e4}
\frac cs g(\Delta u(t,p_0), e^2)=\{\frac {c^2}{s^2}\partial_r u_\theta +\frac cs \partial_r^2 u_\theta\}|_{(t,p_0)}.
\ee
Then \eqref{e2}, \eqref{e3}, \eqref{e4} and \eqref{Deltae3} give
 \begin{align*}
&{\partial_r}|_{p_0} g(u_0,e_2)- {\partial_r}|_{p_0} g(u(t_0, \cdot),e_2)\\
&\quad=\int^{t_0}_{0} \{ \frac 2{s^2}\partial_r u_\theta-\frac cs \partial^2_r u_\theta -\partial_r^3 u_\theta+\frac 1s \partial^2_r\partial_\theta u_r\\
 &\quad\qquad +\frac 1s\partial_\theta \partial^2_r u_r -\frac {c^2}{s^2}\partial_r u_\theta -\frac cs \partial_r^2 u_\theta   \mp2a^2{\partial_r}u_\theta\} \dd t\\
&\quad=\int^{t_0}_{0} (\frac {2-c^2}{s^2}\mp2a^2)\partial_r u_\theta-2\frac cs \partial^2_r u_\theta -\partial_r^3 u_\theta+\frac 2s \partial^2_r\partial_\theta u_r\dd t \\
 &\quad=\int^{t_0}_{0} \frac{c^2}{s^2}\partial_r u_\theta-2\frac cs \partial^2_r u_\theta -\partial_r^3 u_\theta+  \frac 2s\partial_\theta \partial^2_r u_r  \dd t. 
\end{align*}
 Finally, by using \eqref{df} again, we can rewrite the last term as
 \[
 \frac 2s \partial_\theta \partial^2_r u_r|_{\partial K}=-\frac 2{s^2}\partial_r\partial_\theta^2 u_\theta|_{\partial K}.
 \]
 It follows
  \begin{align}\label{final}
& {\partial_r}|_{p_0} g(u_0,e_2)-{\partial_r}|_{p_0} g(u(t_0, \cdot),e_2)\nonumber\\
& =\int^{t_0}_{0} \frac{c^2(\delta)}{s^2(\delta)}\partial_r u_\theta(t,p_0)-2\frac {c(\delta)}{s(\delta)} \partial^2_r u_\theta(t,p_0) -\partial_r^3 u_\theta(t,p_0)-\frac 2{s^2(\delta)}\partial_r\partial_\theta^2 u_\theta(t,p_0)  \dd t ,
\end{align}
or equivalently
\begin{align}\label{finaltoo}
\alpha_1(0)-\alpha_1(t_0)\nonumber=\int^{t_0}_{0}k^2\alpha_1(t)-\alpha_3(t)-2k\alpha_2(t)-2\eta(t) \dd t ,
\end{align}
which gives the desired ODE. 
\section{With Coriolis force and in-flow condition case}\label{Earth}
Recall here we have 
\be\label{inflow}
u_\theta|_{\partial K}=0, \quad u_r|_{\partial K}=\lambda_0\in \R,\quad \partial_r u_r|_{\partial K}=0.
\ee
and we work with \eqref{ceq1}-\eqref{ceq2}.  Then note that 
\begin{equation*}
g(\beta\cos(ar)* u,e_2)=\beta\cos(ar)g(*u,e_2)=\beta\cos(ar)u_r
\end{equation*}
and 
\begin{equation}\label{Coriolis}
\partial_r|_{p_0} g(\beta\cos(ar)* u,e_2)=-a\beta\lambda_0\sin(a\delta).
\end{equation}
Next, recall \eqref{convection1} 
\begin{align*}
g(\nabla_u u, e_2)=u_r\partial_ru_\theta+u_\theta u_r \frac cs+\frac 1s u_\theta\partial_\theta u_\theta.
\end{align*}
Then differentiate, evaluate on the boundary, and this time use \eqref{inflow} to obtain
\[
\partial_r g(\nabla_u u, e_2)|_{\partial K}=\lambda_0\partial^2_r u_\theta  +\lambda_0\partial_ru_\theta\frac cs+ \frac 1s \partial_r u_\theta \partial_\theta u_\theta.
\]
From the divergence free condition and \eqref{inflow} we have 
\be\label{convection3}
\partial_r g(\nabla_u u, e_2)|_{\partial K}=\{\lambda_0\partial^2_r u_\theta  +\lambda_0\partial_ru_\theta\frac cs -\frac 1s \partial_r u_\theta \partial_r(su_r)\}|_{\partial K}=\lambda_0\partial^2_r u_\theta|_{\partial K}.
\ee
Another place where we obtain an extra term is in \eqref{pressureb}, where due to \eqref{convect}, \eqref{inflow} and Coriolis term in \eqref{ceq1}, now we have
\[
\dd p=\Delta u -\lambda_0\partial_r u_\theta e_2+2a^2\lambda_0e_1-\beta\cos(a\delta)\lambda_0e_2
\]
on $\partial K$.
With the Coriolis term \eqref{Coriolis},
then \eqref{e2} becomes
\begin{align}
& {\partial_r}|_{p_0} g(u_0,e_2)-{\partial_r}|_{p_0} g(u(t_0, \cdot),e_2)\nonumber\\
&=-\int^{t_0}_{0}   \{{\partial_r}g(\Delta  u(t, p_0), e^2)-\frac 1s \partial_\theta g(\Delta u(t, p_0), e^1) \nonumber\\
&\qquad\qquad +\frac cs g(\Delta u(t,p_0)-\lambda_0(\partial_r u_\theta+\beta\cos(a\delta))e^2 , e^2) \nonumber\\
&\qquad\qquad+2a^2{\partial_r}u_\theta-\lambda_0\partial^2_r u_\theta
+a\beta\lambda_0\sin(a\delta)\}\dd t.\label{e8}
\end{align}
We then repeat the same computations that followed \eqref{e2}.  The results are the same except that we have the four extra terms that appeared in \eqref{e8}.  This turns \eqref{final} into
  \begin{align*}\label{final2}
& {\partial_r}|_{p_0} g(u_0,e_2)-{\partial_r}|_{p_0} g(u(t_0, \cdot),e_2)\nonumber\\
& =\int^{t_0}_{0}k^2\alpha_1(t)-\alpha_3(t)-2k\alpha_2(t)-2\eta(t) +\lambda_0(k\alpha_1(t)+\alpha_2(t))
-\lambda_0\beta(a\sin(a\delta)-k\cos(a\delta)) \dd t. 
\end{align*}
\noindent
{\bf Acknowledgments.}\
The first author is partially supported by a grant from the National Science Council of Taiwan (NSC 101-2115-M-009-016-MY2).  The second author is partially supported by a grant from the Simons Foundation \#246255.  The third author is partially supported by JSPS KAKENHI Grant Number 25870004.  We also would like to thank our three institutions for the hospitality during the visits when this work was carried out.

\bibliographystyle{amsart}


\end{document}